\documentclass[twoside,12pt]{amsart} 
\usepackage{amsfonts,amssymb,latexsym,amsthm, epsfig,amsmath,oldgerm,amscd}
\def\Com{\UseComputerModernTips}
\renewcommand{\baselinestretch}{1}
\newcommand{\by}[1]{\stackrel{#1}{\rightarrow}}
\newcommand{\longby}[1]{\stackrel{#1}{\longrightarrow}}
\newcommand{\mf}{\mathfrak} 
\newcommand{\surj}{\twoheadrightarrow}
\newcommand{\inj}{\hookrightarrow} 
\newcommand{\lra}{\longrightarrow}
\newcommand{\equi}{\Leftrightarrow}
\newcommand{\ooplus}{\displaystyle{\mathop\oplus}}
\newcommand{\ol}{\overline} 
\newcommand{\wt}{\widetilde}
\newcommand{\iso}{\by \sim} 
\newcommand{\ra}{\rightarrow}
\newcommand{\Ra}{\Rightarrow}
\newcommand{\la}{\leftarrow}
\newcommand{\La}{\leftarrow} 
\newcommand{\eps}{\varepsilon}
\newcommand{\stk}{\stackrel}
\newcommand{\mc}{\mathcal}
\newcommand{\mbb}{\mathbb}
\newcommand{\tn}{\textnormal}
\newcommand{\ubr}{\underbrace}
\newcommand{\obr}{\overbrace}
\newcommand{\myoplus}{\operatornamewithlimits{\oplus}}
\newtheorem{de}{Definition}[section]
\newtheorem{re}[de]{Remark}
 
\newtheorem{tr}[de]{Theorem}
\newtheorem{lm}[de]{Lemma} 

\newtheorem{nt}[de]{Notation} 

\newtheorem{co}[de]{Corollary}

\newcommand{\LMD}{\Lambda}
\newcommand{\lmd}{\lambda}
\newcommand{\lb}{\linebreak}
\def\vp{\rm \vspace{0.2cm}}
\def\M{\rm Max}
\def\m{\mf{m}}
\def\st{\rm St}
\def\sr{\rm sr}
\def\hb{\hfill$\Box$}
\def\aut{\rm Aut} 
\def\sp{\rm Spec} 
\def\tran{\rm Trans}
\def\fl{\rm ETrans}
\def\Um{\rm Um}
\def\GL{\rm GL}
\def\GH{\rm GH}
\def\GQ{\rm GQ}
\def\SL{\rm SL}
\def\EQ{\rm EQ}
\def\EO{\rm EO}
\def\EH{\rm EH}
\def\SH{\rm SH}
\def\SQ{\rm SQ}
\def\SO{\rm SO}
\def\E{\rm E}
\def\EH{\rm EH}
\def\T{\rm T}
\def\ET{\rm ET}
\def\G{\rm G}
\def\Sp{\rm Sp}
\def\ko{\rm K_1O}
\def\ksp{\rm K_1Sp}
\def\k{\rm K_1}
\def\K{\rm K}
\def\KH{\rm KH}
\def\KQ{\rm KQ}
\def\I{\rm I}
\def\T{\rm T}
\def\O{\rm O}
\def\ESp{\rm ESp}
\def\EO{\rm EO}
\def\C{\rm C}
\def\es{\rm S}
\def\sdim{\rm sdim}

\begin{document}
\title{A note on Quadratic and Hermitian Groups}
\author{Rabeya Basu}

\date{}
\maketitle
{\tiny \begin{center}{\it 2000 Mathematics Subject Classification:
{13C10, 11E57, 11E70 {\scshape 19Bxx}}}
\end{center}}
{\tiny \begin{center}{\it Key words: Bilinear forms, Quadratic forms, Hermitian forms}
\end{center}}
{\tiny Abstract: In this article we deduce an analogue of Quillen's 
Local-Global Principle for the elementary subgroup of the general quadratic group and 
the hermitian group. We show that the unstable ${\k}$-groups 
of the hermitian groups are nilpotent by abelian. This generalizes earlier results of 
A. Bak, R. Hazrat, N. Vavilov and {\it etal.}.}

\section{\large Introduction} 

The vigorous study of classical groups, more generally algebraic ${\K}$-theory, stimulated in mid sixties in attempt to give a solution of Serre's Problem for projective modules. This
prominent theorem  in commutative algebra states that {\it finitely generated projective modules 
over a polynomial ring over a field are free}. ({\it cf.} Faisceaux Algebriques Coherent, 1955). 
The beautiful book {\it Serre's Problem on Projective Modules} by T.Y. Lam  gives a comprehensive account of the mathematics surrounding Serre's Problem. Later we see analogue of Serre's Problem for other classical groups in the work of H. Bass, A. Suslin, V.I. Kopeiko, R. Parimala and {\it etal.} in \cite{Bass}, \cite{KOP}, \cite{SUSK}, \cite{P2}. 
We discuss the following problems related to Serre's Problem, {\it viz.} normality property of the elementary subgroup of full automorphism group, Quillen's Local-Global Principle, stability for  ${\k}$-functors, and structure of unstable ${\k}$-groups of classical groups and modules.

From the work of H. Bass ({\it cf.}\cite{B}),  J.S. Wilson ({\it cf.}\cite{W}), L.N. Vaserstein ({\it cf.}\cite{vas2}) and {\it  etal.} it is well known that the elementary subgroup ${\E}_n(A)$ of ${\GL}_n(A)$ 
plays a crucial role in the study of {\it Classical ${\K}$-Theory}. In \cite{V}, L.N. Vaserstein proved that for an associative ring $A$, which is finite over its center, ${\E}_n(A)$ is a normal subgroup of 
${\GL}_n(A)$ for $n\ge 3$. Later analogue results for the symplectic and orthogonal groups were done by Kopeiko and Suslin-Kopeiko in \cite{SUS} and \cite{SUSK} (for respective cases). The difficulties of quadratic Serre's Problem for characteristic 2 was first noted by Bass in \cite{Bass}. In fact, in many cases it is difficult to handle classical groups over fields of characteristic 2, rather than classical groups over fields of characteristic $\ne$ 2. 
(For details see \cite{HV1}). In 1969, A. Bak resolved this problem by introducing {\it Form Rings} and {\it form parameter}. We also see some results in this direction in the work of Kle\u\i n, Mikhalev, Vaserstein {\it etal.} in \cite{K1}, \cite{K2}, \cite{V2}. The concept of {\it form parameter} also appears in the work of K. McCrimmon which plays an important role in his classification theory of Jorgan algebras {\it cf.}\cite{M}, for details see (\cite{HO}, footnote pg. 190.) and  \cite{J}. 

In \cite{BRK}, it has been shown that this normality criterion 
is related to the following well known Local-Global Principle 
introduced by A. Suslin to give a matrix
theoretic proof of Serre's conjecture on projective modules.

{\bf Suslin's Local-Global Principle:} 
{\it Let $R$ be a commutative ring with identity and 
$\alpha(X)\in {\GL}_n(R[X])$ with $\alpha(0)={\I}_n$. If 
$\alpha_{\mf{m}}(X)\in {\E}_n(R_{\mf{m}}[X])$ for every maximal ideal 
$\mf{m}\in {\M}(R)$, then $\alpha(X)\in {\E}_n(R[X])$. } \vp 

We see generalization of this Principle for the symplectic group in \cite{KOP}, and for 
the orthogonal group in \cite{SUSK}. 
In \cite{BRK},  we have shown that the question of normality of the
elementary subgroup of the 
general linear group, symplectic and orthogonal groups, is equivalent 
to the above Local-Global Principle where the base ring is associative with identity and 
is finite over its center. In that article we have treated uniformly above three classical groups. 
Motivated by the work of 
A. Bak, R.G. Swan, L.N. Vaserstein and others  in \cite{BBR} the author with 
A. Bak and R.A. Rao has  established an analogue 
of Suslin's Local-Global Principle for the transvection subgroup of 
the automorphism group of projective, 
symplectic and orthogonal modules of global rank at least 1 and 
local rank at least 3, under the assumption that the projective module has 
constant local rank and that the symplectic and orthogonal modules are locally 
an orthogonal sum of a constant number of hyperbolic planes.
In this article we deduce an analogue Local-Global Principle 
for  the general quadratic (or unitary) and hermitian groups.
We treat these two groups uniformly and give explicit proofs of those results.
We have overcome many technical difficulties which comes in the 
hermitian case due to the elements $a_1,\ldots,a_r$ 
(with respect to these elements we define the hermitian groups). 
We assume $a_1=0$. The rigorous study of hermitian groups can be found in \cite{T}. 

The study of stability for ${\k}$-functors started in mid sixties and first appears in the work of 
Bass-Milnor-Serre. Then this problem was thoroughly studied by L.N. Vaserstein for the symplectic, orthogonal and unitary groups. For details {\it cf.} \cite{V}, \cite{V2}, and \cite{V3}. After almost thirty years of Vaserstein's results in 1998 this problem was revisited for the linear groups over an affine algebra by R.A. Rao and W. van der kallen in \cite{RV}. Finally, the result settled for the general quadratic and hermitian groups by A. Bak, G. Tang and V. Petrov in \cite{Bak2} and \cite{Bak3}. It has been observed that 
over a regular affine algebra Vaserstein's bounds for the stabilization can be improved for the 
transvection subgroup of full automorphism group of projective, and 
symplectic modules. But cannot be improved for the orthogonal case in general. For details {\it cf.}
\cite{BR}, \cite{BRS}. 

Though the study of stability for ${\k}$-functors started in mid sixties, the structure of unstable ${\k}$-group first studied by A. Bak in 1991 ({\it cf.}\cite{Bak}). He showed that 
the group ${\GL}_n(R)/{\E}_n(R)$ is nilpotent-by-abelian for $n\ge 3$. 
In \cite{HV}, R. Hazrat and N. Vavilov generalized his result for Chevalley groups with irreducible root system. They have shown: Let $\Phi$ be a reduced
irreducible root system of rank $\geq 2$ and $R$ be a commutative ring
such that its Bass-Serre dimension $\delta(R)$ is finite. Then for any   
Chevalley group ${\G}(\Phi, R)$ of type $\Phi$ over $R$ the quotient 
${\G}(\Phi, R)/{\E}(\Phi, R)$ is nilpotent-by-abelian. In particular, 
${\k}(\Phi, R)$ is nilpotent of class at most $\delta(R) + 1$. 
They use the localization-completion method of A. Bak in \cite{Bak}. 
In \cite{BBR}, the author with Bak and Rao give a uniform proof for 
the transvection subgroup of 
full automorphism group of projective, 
symplectic and orthogonal modules of global rank at least 1 and 
local rank at least 3. Our method of proof shows that for classical groups the localization part suffices.
Recently, in ({\it cf.}\cite{BHV}) Bak, Vavilov and Hazrat proved the relative case 
for the unitary and Chevalley groups. But, in my best knowledge, so far there is no result for hermitian groups. I observe that using the above Local-Global Principle,
arguing as in \cite{BBR}, it follows that the unstable ${\k}$ of hermitian group is nilpotent-by-abelian. We follow the line of Theorem 4.1 in \cite{BBR}. 

The first part of the paper serve as an introduction. In section 3 and 4, we discuss Suslin's Local-Global Principle and its equivalence with the normality property of the elementary subgroup of full automorphism group. Finally, in section 4 we study the nilpotent property of unstable ${\KH}_1$. 

\section{Preliminaries} Let $R$ be an associative ring with identity.
Recall that an involutive anti-homomorphism
({\bf involution}, in short) is a homomorphism $*:R\ra R$ such that  
$(x-y)^{*}=x^{*}-y^{*}$, $(xy)^{*}=y^{*}x^{*}$ and $(x^{*})^{*}=x$ for any 
$x,y\in R$. For any left $R$-module $M$ the involution induces a left module 
structure to the right $R$-module $M^{*}$=Hom$(M,R)$ given by 
$(xf)v=(fv)x^{*}$, where $v\in M$, $x\in R$ and $f\in M^{*}$.
Let $-:R\ra R$ be an involution for which there is an element $\lmd\in R$ such 
that $\lmd a \ol{\lmd}=\ol{\ol{a}}$ for all $a\in R$. Setting $a=1$ we 
obtain $\lmd \ol{\lmd}=1$. Direct computation shows that $\ol{\lmd}\lmd=1$, 
from which it will follow that $\lmd$ is invertible in $R$ and 
$\lmd^{-1}=\ol{\lmd}$. The element $\lmd$ is called a symmetry of the
involution ``$-$'' and it is unique up to an element $c\in C(R)$, where $C(R)$
is the center of $R$, such that $c\ol{c}=1$. An involution with symmetry
$\lmd$ will be called a $\lmd$-involution.

We assume that $-:R\ra R$ denote a $\lmd$-involution on $R$. Let 
$\tn{max}^{\lmd}(R)=\{a\in R\,|\, a=-\lmd\ol{a}\}$ and $\tn{min}^{\lmd}(R)=\{
a-\lmd\ol{a}\,|\, a\in R\}$. One checks that $\tn{max}^{\lmd}(R)$ and 
$\tn{min}^{\lmd}(R)$ are closed under the operation $a\mapsto \ol{x}ax$ for any
$x\in R$. A $\lmd$-form parameter on $R$ is an additive subgroup $\LMD$ of $R$
such that $\tn{min}^{\lmd}(R)\subseteq \LMD\subseteq\tn{max}^{\lmd}(R)$, and 
$\ol{x}\LMD x\subseteq \LMD$ for all $x\in R$. Note that if $\lambda=-1$, then 
$\tn{max}^{\lmd}(R)=R$. The symplectic case is when $\LMD=R$. On the other hand if 
$\lambda=1$, then $\tn{min}^{\lmd}(R)=0$. The orthogonal case is when $\LMD=0$. 

Let $R$ possesses a $\lmd$-involution 
$-:a\mapsto \ol{a}$, for $a\in R$. For a matrix $M=(m_{ij})$ over $R$ we  
define $\ol{M}=(\ol{m}_{ij})^t$. Let $A_1$ denotes the diagonal matrix 
$[a_1,\ldots,a_r]$ for $a_1,\ldots,a_r\in R$, and $A=A_1\perp I_{n-r}$. 
We define the forms 
{\tiny $$\psi^q_n= \begin{pmatrix} 0 & \lmd I_n \\ I_n &0\end{pmatrix} \,\,\, \tn{  and    }  
\psi^h_n= \begin{pmatrix} A & \lmd I_n \\ I_n &0\end{pmatrix}.$$ }

\begin{de} {\bf General Quadratic Group
${\GQ}_{2n}(R,\LMD)$:} \tn{The group generated by the all non-singular 
$2n\times 2n$ matrices $\{\sigma\in {\GL}_{2n}(R)\,|\, \ol{\sigma}\psi^q_n
\sigma=\psi^q_n\}$, where $\sigma=\begin{pmatrix} \alpha & \beta \\ \gamma &
\delta \end{pmatrix}$ ($\alpha, \beta, \gamma, \delta$ are $n\times n$ block 
matrices) and $\ol{\gamma} \alpha,\ol{\delta}\beta\in \LMD$.}
\end{de} 
\begin{de} {\bf General Hermitian Group of the elements $a_1,\!\!\ldots \!\!,a_r$\lb
${\GH}_{2n}(R, a_1,\ldots ,a_r, \LMD)$:} \tn{The group generated by the all 
non-singular $2n\times 2n$ matrices $\{\sigma\in 
{\GL}_{2n}(R)\,|\, \ol{\sigma}\psi^h_n \sigma=\psi^h_n\}$.}
\end{de} 

A typical element in ${\GQ}_{2n}(R,\LMD)$ and ${\GH}_{2n}(R, a_1,\ldots,a_r, \LMD)$ 
is denoted by a $2n\times 2n$ matrix {\tiny $\begin{pmatrix} \alpha & \beta \\ 
\gamma &\delta \end{pmatrix}$} where $\alpha, \beta, \gamma, \delta$ are
$n\times n$ block matrices. There is a standard embedding, ${\GQ}_{2n}(R,\LMD)\ra {\GQ}_{2n+2}(R,\LMD)$ given by 
{\tiny $$\begin{pmatrix} \alpha & \beta \\ 
\gamma &\delta \end{pmatrix} \mapsto \begin{pmatrix}\alpha & 0 & \beta & 0\\
0 & 1 & 0 & 0\\ \gamma & 0 & \delta & 0\\0 & 0 & 0 & 1\end{pmatrix}$$}
called the stabilization map. This allows us to identify ${\GQ}_{2n}(R,\LMD)$ with a subgroup in ${\GQ}_{2n+2}(R,\LMD)$. Similarly, there is an obvious embedding
${\GH}_{2n}(R, a_1,\ldots,a_r, \LMD)\ra {\GH}_{2n+2}(R, a_1,\ldots,a_r, \LMD)$. 

We recall the definition of the {\bf elementary quadratic matrices} in 
${\GQ}_{2n}(R,\LMD)$. 
Let $\rho$ be the permutation defined by $\rho(i)=n+i$ for $i=1,\ldots,n$. 
Let $e_i$ denote the vector with $1$ in the $i$-th position and O's elsewhere. 
Let $e_{ij}$ be the matrix with $1$ in the $ij$-th position and 0's
elsewhere. We define 
{\tiny $$q\eps_{ij}(a)=I_{2n}+ae_{ij}-\ol{a}e_{\rho(j)\rho(i)} \,\,\tn{ for }
\,\,a\in R, \,\,\tn{ and } 1\le i, j\le n, i\ne j,$$
$$qr_{ij}(a)=I_{2n}+ae_{i\rho(j)}-\lmd\ol{a}e_{j\rho(i)} \,\,\tn{ for }
\,\,a\in R, \,\,\tn{ and } 1\le i, j\le n,$$ }
(Note that, if $i=j$, it forces that $a\in \LMD$), and  
{\tiny $$ql_{ij}(a)=I_{2n}+ae_{\rho(i)j}-\ol{\lmd}\ol{a}e_{\rho(j)i} \,\,\tn{ for }
\,\,a\in R, \,\,\tn{ and } 1\le i, j\le n.$$}
(Note that, if $i=j$, it forces that $a\in \ol{\LMD}$). One checks that these
above matrices belong to ${\GQ}_{2n}(R,\LMD)$ when $a\in R$. 
\begin{de} {\bf n'th Elementary Quadratic Group
${\EQ}_{2n}(R,\LMD)$:} \tn{The subgroup generated by $q\eps_{ij}(a)$, 
$qr_{ij}(a)$ and $ql_{ij}(a)$ for $a\in R$ and $1\le i,j\le n$.} 
\end{de} 
To define {\bf elementary Hermitian matrices}, we need to consider the set 
$C=\{(x_1,\ldots,x_r)^t\in (R^r)^t\,|\,\underset{i=1}{\overset{r}\sum}
\ol{x}_ia_ix_i\in \tn{min}^{-\lmd}(R)\}$
for $a_1,\ldots,a_r$ as above. In order to overcome the technical difficulties 
caused by the elements $a_1,\ldots,a_r$, we shall finely partition a typical
matrix {\tiny $\begin{pmatrix} \alpha & \beta \\ 
\gamma &\delta \end{pmatrix}$ of ${\GH}_{2n}(R, a_1,\ldots,a_r, \LMD)$} into the form 
{\tiny $$\begin{pmatrix} \alpha_{11} & \alpha_{12} & \beta_{11} & \beta_{12} \\
\alpha_{21} & \alpha_{22} & \beta_{21} & \beta_{22}\\
\gamma_{11} & \gamma_{12} & \delta_{11} & \delta_{12} \\
\gamma_{21} & \gamma_{22} & \delta_{21} & \delta_{22} \end{pmatrix}$$  }
where $\alpha_{11},\beta_{11}, \gamma_{11}, \delta_{11}$ are $r\times r$
matrices, $\alpha_{12},\beta_{12},\gamma_{12},\delta_{12}$ are $r\times
(n-r)$ matrices, $\alpha_{21},\beta_{21},\gamma_{21},\delta_{21}$
are $(n-r)\times r$ matrices, and $\alpha_{22},\beta_{22},\gamma_{22},
\delta_{22}$ are $(n-r)\times (n-r)$ matrices. By (\cite{T}, Lemma 3.4), 
\begin{eqnarray} \label{tn1} 
\tn{the columns of}
\alpha_{11}-I_r,\alpha_{12},\beta_{11},\beta_{12},\ol{\beta}_{11},
\ol{\beta}_{21}, \ol{\delta}_{11}-I_r,\ol{\delta}_{21} \in C. 
\end{eqnarray} 
One checks straightforward that the subgroup of ${\GH}_{2n}(R, a_1,\ldots,a_r, \LMD)$ 
consisting of 
{\tiny \begin{eqnarray*} \label{tn2}
 \left\{ \begin{pmatrix} I_r & 0 & 0 & 0\\
0 & \alpha_{22} & 0 & \beta_{22}\\
0 & 0 & I_r & 0 \\
0 & \gamma_{22} & 0 & \delta_{22} \end{pmatrix}
\in {\GH}_{2n}(R, a_1,\ldots,a_r) \right\}\cong {\GH}_{2(n-r)}(R,a_1,\!\!\ldots\!\!,a_r,\LMD).
\end{eqnarray*}}
The first three kinds of generators are taken for the most part from 
${\GQ}_{2(n-r)}(R,\LMD)$, which is embedded as above as a subgroup of
${\GH}_{2n}$ and the last two kinds are motivated by the result \eqref{tn1}
concerning the column of a matrix in ${\GH}_{2n}$. We define 
{\tiny $$h\eps_{ij}(a)=I_{2n}+ae_{ij}-\ol{a}e_{\rho(j)\rho(i)} \,\,\tn{for}
\,\,a\in R, \tn{ and } r+1\le i\le n, 1\le j\le n, i\ne j,$$
$$hr_{ij}(a)=I_{2n}+ae_{i\rho(j)}-\lmd\ol{a}e_{j\rho(i)} \,\,\tn{for}
\,\,a\in R, \tn{ and } r+1\le i, j\le n,$$}
(Note that, if $i=j$, it forces that $a\in \LMD$) and 
{\tiny $$hl_{ij}(a)=I_{2n}+ae_{\rho(i)j}-\ol{\lmd}\ol{a}e_{\rho(j)i} \,\,\tn{for}
\,\,a\in R, \tn { and } 1\le i, j\le n.$$}
(Note that, if $i=j$, it forces that $a\in \ol{\LMD}$). It follows from 
\eqref{tn2} that these matrices belong to ${\GH}_{2n}(R,a_1,\ldots,a_r, \LMD)$ 
when $a\in R$. 

For $\zeta\in (x_1,\ldots,x_r)^t\in C$, let $\zeta_f\in R$ such that 
$\zeta_f+\lmd\ol{\zeta}_f=\underset{i=1}{\overset{r}\sum}
\ol{x}_ia_ix_i$. (The element $\zeta_f$ is not unique in general). We define
{\tiny $$hm_i(\zeta)=\begin{pmatrix} I_r & \alpha_{12} & 0 & 0\\
0 & I_{n-r} & 0 & 0\\
0 & -\ol{A}_1\alpha_{12} & I_r & 0 \\
0 & \gamma_{22} & -\ol{\alpha}_{12} & I_{n-r} \end{pmatrix}\,\,\tn{for} 
\,\,(\zeta\in C \,\,\& \,\, r+1\le i\le n)$$}
a $2n\times 2n$ matrix, where $\alpha_{12}$ is the $r\times (n-r)$ matrix with
$hm_i(\zeta)e_{i-r}=\zeta$ and all other column's zero, and $\gamma_{22}$ is
the $(n-r)\times (n-r)$ matrix with $\gamma_{22}e_{i-r,\,\,i-r}=\ol{\zeta}_f$ 
and $0$ elsewhere. 

As above, we define 
{\tiny $$hr_i(\zeta)=\begin{pmatrix} I_r & 0 & 0 & \beta_{12}\\
0 & I_{n-r} & -\lmd\ol{\beta}_{12} & \beta_{22}\\
0 & 0 & I_r & -\ol{A}_1\beta_{12} \\
0 & 0 & 0 & I_{n-r} \end{pmatrix}\,\,\tn{ for } 
\,\,(\zeta\in C \,\,\& \,\, r+1\le i\le n)$$}
a $2n\times 2n$ matrix, where $\beta_{12}$ is the $r\times (n-r)$ matrix with
$hr_i(\zeta)e_{i-r}=\zeta$ and all other column's zero, and $\beta_{22}$ is
the $(n-r)\times (n-r)$ matrix with $\beta_{22}e_{i-r,\,\,i-r}=\ol{\zeta}_f$ 
and $0$ elsewhere. 
\begin{de} \tn{Each of the above matrices is called an 
{\bf elementary hermitian matrix} for the elements $a_1,\ldots a_r$. }
\end{de}

{\tiny Note that if $\eta=e_{pq}(a)$ is an elementary generator in ${\GL}_s(R)$, 
then the matrix $(I_{n-s}\perp
\eta\perp I_{n-s} \perp \eta^{-1})\in h\eps_{ij}(a)$. 
It has been shown in (\cite{T}, $\S 5$) that each of the 
above matrices is in ${\GH}_{2n}(R, a_1,\ldots,a_r, \LMD)$.}
\begin{de} {\small {\bf n'th Elementary Hermitian Group 
{\tiny${\EH}_{2n}(R,a_1,\ldots,a_r, \LMD)$:}}} \tn{The group generated by $h\eps_{ij}(a)$, 
$hr_{ij}(a)$, $hl_{ij}(a)$, $hm_i(\zeta)$ and $hr_i(\zeta)$ for $a\in R$ 
and $1\le i,j\le n$.} 
\end{de} 
{\tiny {\bf Blanket Assumption:} We assume that $2n\ge 6$ and  
$n>r$ while dealing with the hermitian case. We consider only isotropic
vectors of $\LMD$. We do not want to put any restriction on the
elements of $C$. Therefore we assume that $a_i\in \tn{min}^{\lmd}(R)$ for 
$i=1,\ldots,r$, as in that case $C=R^r$. When dealing with the Hermitian case we always assume 
$a_1=0$.}
\begin{nt} \tn{In the sequel \tn{M}$(2n,R)$ will denote the set of all 
$2n\times 2n$ matrices. $\G(2n,R, \LMD)$ will denote either the quadratic group 
${\GQ}_{2n}(R,\LMD)$ or the hermitian group ${\GH}_{2n}(R,a_1,\ldots,a_r, \LMD)$.
${\es}(2n, R, \LMD)$ will denote respective subgroups ${\SQ}_{2n}(R,\LMD)$ or ${\SH}_{2n}(R,a_1,\ldots,a_r,\LMD)$ with matrices of determinant 1, in the case when $R$ will be commutative.
${\E}(2n,R, \LMD)$ will denote the corresponding elementary subgroups
${\EQ}_{2n}(R,\LMD)$ and ${\EH}_{2n}(R,a_1,\ldots,a_r, \LMD)$. 
To treat uniformly we denote the generators of ${\E}(2n,R,\LMD)$ by $ge_{ij}(a)$ for $a\in R$.
Let $\LMD[X]$ denote the $\lambda$-form parameter on $R[X]$ induced from $(R,\LMD)$, and
let $\LMD_s$ denote the $\lambda$-form parameter on $R_s$ induced from $(R,\LMD)$. }
\end{nt} 

For any column vector $v\in R^{2n}$ we denote
$\widetilde{v}_q=\ol{v}^t\psi^q_n$ and $\widetilde{v}_h=\ol{v}^t\psi^h_n$.
\begin{de} \tn{We define the map $M:R^{2n} \times R^{2n} \ra \tn{M}(2n,R)$ and 
the inner product $\langle \,,\rangle$ as follows: }
{\tiny \begin{align*}
M(v,w) & = v.\widetilde{w}_q-\ol{\lambda}\ol{w}.\widetilde{v}_q \,\,\,\, \tn{ when } 
      {\G}(2n,R)={\GQ}_{2n}(R,\LMD)\\
  & = v.\widetilde{w}_h -\ol{\lambda}\ol{w}.\widetilde{v}_h \,\,\,\,  \tn{ when } 
      {\G}(2n,R)={\GH}_{2n}(R,a_1,\ldots,a_r, \LMD),\\ 
\langle v,w\rangle & = \widetilde{v}_q.w \,\,\,\, \tn{ when } 
{\G}(2n,R)={\GQ}_{2n}(R,\LMD) \\
  & = \widetilde{v}_h.w \,\,\,\, \tn{ when } 
             {\G}(2n,R)={\GH}_{2n}(R,a_1, \ldots , a_r, \LMD).
\end{align*} }
\end{de} 

We recall following very well known facts: 

\begin{lm} \label{nor} $($cf. \cite{Bak1}, \cite{T}$)$
The group  ${\E}(2n,R, \LMD)$ is perfect for $n\ge 3$ in the quadratic case and for $n\ge r+3$ in the hermitian case, {\it i.e.} 
{\tiny $$[{\E}(2n,R, \LMD), {\E}(2n,R, \LMD)] = {\E}(2n,R, \LMD).$$ }
\end{lm}

\begin{lm} {\bf (Splitting property):} For all $x, y \in R$ 
{\tiny $$g_{ij}(x+y)=g_{ij}(x)g_{ij}(y).$$} 
\end{lm}

\begin{lm}
Let $G$ be a group, and $a_i$, $b_i \in G$, for 
$i = 1, \ldots, n$. Then  for $r_i={\underset{j=1}
{\overset{i}\Pi}} a_j$ we have
${\underset{i=1}{\overset{n}\Pi}}a_i b_i=
{\underset{i=1}{\overset{n}\Pi}}r_ib_ir_i^{-1} 
{\underset{i=1}{\overset{n}\Pi}}a_i.$
\end{lm}
\begin{nt} \tn{ By ${\G}(2n,R[X], \LMD [X], (X))$ we mean the group of all invertible matrices
over $R[X]$ which are ${\I}_n$ modulo $(X)$.} 
\end{nt} 
\begin{lm} \label{key1} The group 
${\G}(2n,R[X],\LMD [X], (X)) \cap {\E}(2n,R[X], \LMD [X])$
is generated by the elements of the type
$ \eps ge_{ij}(Xh(X))\eps^{-1}$, where  $\eps \in {\E}(2n,R, \LMD)$,
$h(X)\in R[X]$.
\end{lm}
\section{Suslin's Local-Global Principle} 

In his remarkable thesis ({\it cf.}\cite{Bak1}) Bak showed that for 
a form ring $(R, \LMD)$ the elementary subgroup 
${\EQ}(2n,R, \LMD)$ is perfect for $n\ge 3$ and hence a normal subgroup of ${\GQ}(2n,R, \LMD)$.
As we have noted earlier that this question is related to 
Suslin's Local-Global Principle for 
the elementary subgroup. In \cite{T}, Tang has shown that for $n\ge r+3$ the elementary hermitian group ${\EH}_{2n}(R,a_1, \ldots , a_r, \LMD)$ is perfect for $n\ge r+3$ and hence a normal subgroup of ${\GH}_{2n}(R,a_1, \ldots , a_r, \LMD)$.
In this section we deduce 
analogue Local-Global Principle for the elementary subgroup of the general quadratic and hermitian groups when $R$ is finite over its center.
We use this result in $\S 5$ to prove the nilpotent property of unstable ${\KH}_1$ 
Furthermore, we show that
if $R$ is finite over its center the normality of elementary subgroup is
equivalent to Local-Global Principle. This generalizes our result in \cite{BRK}. 

The following is the key Lemma, and it tells us the reason why we need to assume that the size of the matrix is at least 6. In \cite{BRK} proof is given for the same result for the linear group. Arguing in similar manner by using identities of commutator laws result follows in the unitary and hermitian cases. A list of commutator laws for elementary generators is stated in (\cite{Bak1}, \S, pg. 43-44) for the unitary group and in (\cite{T}, pg. 237-239) for the hermitian group.

\begin{lm} \label{key3}
 For $m>0$, there are $h_t(X,Y,Z)\in R[X,Y,Z]$ such 
that
{\tiny $$ge_{pq}(Z)ge_{ij}(X^{2m}Y)ge_{pq}(-Z)= \underset{t=1}{\overset{k}\Pi}
ge_{p_tq_t}(X^m h_t(X,Y,Z)).$$ }
\end{lm}

\begin{co} \label{key4}
 If $\eps=\eps_1\eps_2\cdots\eps_r$, where each $\eps_j$ is an
elementary generator, then there are $h_t(X,Y)\in R[X,Y]$ such that
{\tiny $$\eps ge_{pq}(X^{2^rm}Y)\eps^{-1}=\underset{t=1}{\overset{k}\Pi}
ge_{p_tq_t}(X^m h_t(X,Y)).$$}
\end{co}

\begin{lm} \label{key5} 
Let $(R, \LMD)$ be a form ring and $v\in {\E}(2n,R, \Lambda)e_{2n}$. Let $w\in R^{2n}$ be a 
column vector such that $\langle v,w\rangle=0$. Then
 $I_{2n}+M(v,w)\in {\E}(2n,R, \lambda)$.
\end{lm}
{\bf Proof.} 
Let $v=\eps e_{2n}$, where $\eps \in {\E}(2n,R,\LMD)$. Then it follows 
that
$I_{2n}+M(v,w)=\eps (I_{2n}+ M(e_{2n},w_1))\eps^{-1},$ where 
$w_1 = \eps^{-1}w$.
Since $\langle e_{2n},w_1\rangle=\langle v,w\rangle=0$, we get
$w_1^t=(w_{11},\dots,w_{1 \,n-1},0,w_{1\,n+1},\ldots,w_{1 \,2n})$.
Therefore, as $\lmd\ol{\lmd}=\ol{\lmd}\lmd=1$, 
{\tiny $$I_{2n}+M(v,w)= \begin{cases} 
\underset{1\le i\le n-1}
{\underset{1\le j\le n} \Pi}
\,\eps \,ql_{in}(-\ol{\lmd} \ol{w}_{n+1\, i})\,q\eps_{jn}(-\ol{\lmd}\ol{w}_{1\, j})ql_{n\, n}^{-1}(*)
{\eps}^{-1} \\
\underset{1\le i\le n-1}
{\underset{r+1\le j\le n}{\underset{1\le k\le r} \Pi}}
\eps hl_{in}(-\ol{\lmd} \ol{w}_{n+1\, i})\,h\eps_{jn}(-\ol{\lmd}\ol{w}_{1\, j})hm_n
(-\ol{w}_{1 \,k)})hl_{n\, n}^{-1}(*)
{\eps}\!\!\!^{-1} 
\end{cases}$$ }
(in the  quadratic and hermitian cases respectively),
where $\ol{w}_{1 \,\,n+k}=
(w_{1 \,n+k},0,\ldots,0)$. 
Hence the result follows. \hb \vp 

Note that the above implication is true for any associative ring with identity. 
From now onwords we assume that $R$ is finite over its center $C(R)$. 
\begin{lm} \label{Di}
{\bf (Dilation Lemma)} 
Let $\alpha(X)\in {\G}(2n,R[X], \LMD [X])$, with $\alpha(0)=I_{2n}$. If $\alpha_s(X)\in {\E}(2n,R_s[X], \LMD_s [X])$ for some non-nilpotent $s\in R$, then $\alpha(bX) \in  {\E}(2n,R[X], \LMD [X])$ 
for $b\in s^l {\C}(R)$, $l\gg 0$. 

$($Actually, 
we mean there exists some $\beta(X)\in {\E}(2n,R[X], \LMD [X])$ such that
$\beta(0)=I_{2n}$ and $\beta_s(X)=\alpha(bX)$. But, since there is no 
ambiguity, for simplicity we are using the notation $\alpha(bX)$ instead of 
$\beta_s(X)).$
\end{lm} 
{\bf Proof.} We prove the Lemma in 2 steps. \vp \\
{\bf Step 1.} First we prove the following: \vp 

{\it If $\alpha(X)=I_{2n}+ X^d M(v,w)$, where $v\in {\E}(2n,R, \LMD)e_{2n}$ and
$\langle v,w\rangle=0$, then $\alpha(X)\in {\E}(2n,R[X], \LMD [X])$ and it can be
written as a product 
decomposition of the form $\Pi ge_{ij}(Xh(X))$ for $d\gg 0$.}

Replacing $R$ by $R[X]$ in Lemma \ref{key5} we get that
$I_n+XM(v,w)\in {\E}(n,R[X], \LMD[X])$. Let $v=\eps e_1$, where $\eps \in {\E}(n,R)$. 
Since $\lmd\ol{\lmd}=\ol{\lmd}\lmd=1$, as in the proof of Lemma \ref{key5}, we can write 
{\tiny $$I_{2n}+M(v,w)=\begin{cases} 
\underset{1\le i\le n-1}
{\underset{1\le j\le n} \Pi}
\eps ql_{in}(-\ol{\lmd}X \ol{w}_{n+1\, i})\,q\eps_{jn}(-\ol{\lmd}X\ol{w}_{1\, j})ql_{n\, n}^{-1}(*)
{\eps}^{-1} \\
\underset{1\le i\le n-1}
{\underset{r+1\le j\le n}{\underset{1\le k\le r} \Pi}}
\eps hl_{in}(-\ol{\lmd}X \ol{w}_{n+1\, i})\,h\eps_{jn}(-\ol{\lmd}X\ol{w}_{1\, j})hm_n
(-X\ol{w}_{1 \,k)})hl_{n\, n}^{-1}(*){\eps}\!\!\!^{-1} 
\end{cases}$$ }
(in the quadratic and hermitian cases respectively),
where $\ol{w}_{1 \,\,n+k}=
(w_{1 \,n+k},0,\ldots,0)$. 

Now we split the proof into following two cases. 

Case I: $\eps$ is an elementary generators of the type $ge_{pq}(x)$, 
$x\in R$. First applying the homomorphism 
$X\mapsto X^2$ and then applying Lemma \ref{key3} over $R[X]$ we
get $I_n + X^2M(v,w)=\underset{j}\Pi\left(\underset{t=1}{\overset{k}\Pi}
ge_{p_{j(t)}q_{j(t)}}(Xh_{j(t)}(X))\right)$, where $h_{j(t)}(X)\in R[X]$. 
Now the result follows for $d\gg 0$. 

Case II: $\eps$ is a product of elementary generators of the type 
$ge_{pq}(x)$. Let $\mu(\eps)=r$. First 
applying the homomorphism $X\mapsto X^{2^r}$ 
and then applying the Corollary \ref{key4} we see that 
the result is true for $d\gg 0$. \vp \\
{\bf Step 2.} Let $\alpha_s(X)=\underset{k}\Pi ge_{i(k)j(k)}(h_k(X))$, 
where $h_k(X)\in R_s[X]$. So, 
$\alpha_s(XT^d)=\underset{k}\Pi ge_{i(k) j(k)}(h_k(XT^d))$.
Choose $d\ge 2^r$ and $r=\mu(\alpha(X))$. Since $\alpha(0)=I_n$, 
as in Lemma \ref{key1} we get 
\begin{eqnarray*}
\alpha_s(XT^d)  &=& \underset{k}\Pi 
\eps_k ge_{i(k)j(k)}(XT^d\lambda_k(XT^d))\eps_k^{-1}\\ 
&=&\underset{k}\Pi (I_n + XT^d\lambda_k(XT^d)\eps_k M(e_{i(k)}, e_{\sigma(j)(k)})
\eps_k^{-1}),
\end{eqnarray*}
for $\lambda_k(XT^d) \in R_s[X, T]$, $\eps_k \in {\E}(n, R_s, \LMD_s)$. 
Let $v_k=\eps_k e_{i(k)}$. Then taking 
$w_k(X,T)=\ol{\lambda}X\lambda_k(XT^d)\eps_k e_{\sigma(j)(k)}$ otherwise 
we get  $\langle v_k,w_k(X,T)\rangle=0$ 
(without loss of generality we are assuming $\sigma(i)\ne j$). 
Applying result in Step I over the polynomial ring $(R_s[X])[T]$ we get
$I_n+T^d M(v_k,w_k(X,T)) \lb \in {\E}(n,R_s[X,T], \LMD_s[X,T])$, and can be expressed
as a product of the form $\Pi ge_{p_{k(t)} q_{k(t)}}(T h_k(X,T))$, where 
$h_k(X,T)\in R_s[X,T]$. Let $l$ be the maximum of the powers occurring in
the denominators of $h_k(X,T)$ for all $k$. 
Now applying the homomorphism $T\mapsto s^mT$ for 
$m\ge l$ we get $\alpha(bXT^d)\in {\E}(n,R[X,T], \LMD[X,T])$ for some $b\in (s^l)C(R)$. 
Finally, putting $T=1$ we get the required result. \hb

\begin{tr} \label{LG}
{\bf (Local-Global Principle)} 
If $\alpha(X)\in {\G}(2n,R[X], \LMD[X])$, $\alpha(0)={\I}_n$ and 
$\alpha_{\m}(X)\in {\E}(n,R_{\m}[X], \LMD_{\mf m}[X])$ for every maximal ideal $\m \in \M \, (C(R))$, then $\alpha(X)\in {\E}(2n,R[X], \LMD[X])$. $($Note that $R_{\m}$ denotes $S^{-1}R$, where 
$S = C(R) - \m$.$)$
\end{tr}
{\bf Proof.} Since $\alpha_{\m}(X)\in {\E}(2n,R_{\m}[X], \LMD_{\mf m}[X])$ for all
$\m \in {\M}(C(R))$, for each $\m$ there exists $s\in C(R) - \m$ such that 
$\alpha_s(X)\in  {\E}(n,R_s[X], \LMD_s[X])$. Let $\theta(X,T)=
\alpha_s(X+T)\alpha_s(T)^{-1}$. 
Then {\tiny $$\theta(X,T)\in {\E}(2n,(R_s[T])[X], \LMD_s[T][X])$$} and $\theta(0,T)={\I}_n$. 
Then by Dilation Lemma, applied with base ring $R[T]$, 
{\tiny $$\theta(bX,T)\in {\E}(2n,R[X,T], \LMD[X,T]) \tn{ for some } b\in (s^l)C(R),\, l\gg 0.
\hfill\tn{ \,\, \,\,\,}(A)$$}
Let $b_1,b_2,\dots,b_r\in C(R)$ be such that (A) holds and $b_1+\cdots+b_r=1$. 
Then $\theta(b_iX,T)\in {\E}(2n,R[X,T], \LMD[X,T])$ and hence 
$\underset{i=1}{\overset{r}\Pi}\theta(b_iX,T)\in {\E}(2n,R[X,T], \LMD[X,T])$. But, 
{\tiny $$\alpha(X)=\left(\underset{i=1}{\overset{r-1}\Pi} \theta(b_iX,T)
{\mid}_{T=b_{i+1}X+\cdots +b_rX}\right) \theta(b_rX,0).$$ } Since 
$\alpha(0)={\I}_n$, $\alpha(X)\in {\E}(2n,R[X], \LMD[X])$.  \hb

\section{Equivalence of Normality and Local-Global} 

Next we show that if $R$ is a commutative ring with identity
and $A$ is an associative $R$-algebra such that $A$ is finite 
as a left $R$-module, then the normality criterion of elementary 
subgroup is equivalent to Suslin's Local-Global Principle for 
above two classical groups. 

One of the crucial ingredients in the proof of the above theorem is the 
following result which states that the group ${\E}$ acts transitively on 
unimodular vectors. The precise statement of the fact is the following: 
\begin{tr} \label{swan}
Let $R$ be a semilocal ring $($not necessariy comuttive$)$ with involution and 
$v=(v_1,\ldots,v_{2n})^t$ be a  unimodular and istropic vector in $R^{2n}$.
Then $v\in {\E}(2n, R)e_{2n}$ for $n\ge 2$. {\it i.e.} ${\E}(2n,R)$ acts
transitively on ${\Um}_{2n}(R)$. 
\end{tr}

Let us first recall some known facts before we give a proof of the theorem. 

\begin{de} \tn{An associative ring $R$ is said to be {\bf semilocal} if 
$R/\tn{rad}(R)$ is artinian semisimple.}
\end{de} 

We recall the following three lemmas. . 
\begin{lm} \label{HB}
\tn{(H. Bass)} 
Let $A$ be an associative $R$-algebra such that $A$ is finite as a left 
$R$-module and $R$ be a commutative local ring with identity. 
Then $A$ is semilocal. 
\end{lm} 
{\bf Proof.} Since $R$ is local, $R/\tn{rad} (R)$ is a division ring by 
definition. That implies $A/\tn{rad} (A)$ is a finite module over the division 
ring $R/\tn{rad}(R)$ and hence is a finitely generated vector space. Thus 
$A/\tn{rad} (A)$ artinian as $R/\tn{rad}(R)$ module and hence  
$A/\tn{rad}(A)$ artinian as $A/\tn{rad}(A)$ module, so 
it is an artinian ring.

It is known that an artin ring is semisimple if its radical
 is trivial. Thus $A/\tn{rad}(A)$ is semisimple, as 
$\tn{rad}(A/\tn{rad}(A))=0$. 
Hence $A/\tn{rad}(A)$ artinian semisimple. Therefore, $A$ is semilocal by 
definition. \hb 
\begin{lm} \tn{(H. Bass) (\cite{B}, Lemma 4.3.26)} \label{B}
Let $R$ be a semilocal ring $($may not be commutative$)$,  and let 
$I$ be a left ideal of $R$. Let $a$ in $R$
be such that $Ra+I=R$. Then the coset $a+I=\{a+x \,|\,x\in I\}$
contains a unit of $R$. 
\end{lm}
{\bf Proof.} We give a proof due to R.G. Swan. 
We can factor out the radical and assume that $R$ is semisimple 
artinian. Let $I=(Ra\cap I)\oplus I'$. Replacing $I$ by $I'$ we can assume 
that $R=Ra\oplus I$. Let $f:R\ra Ra$ by $r\mapsto ra$ for $r\in R$. Therefore, 
we get an split exact sequence $0\lra J\lra R \stk{f}\lra Ra\lra 0$, for some 
ideal $J$ in $R$ which gives us a map $g:R\ra J$ such that $R \stk{(f,g)}\lra 
Ra\oplus J$ is an isomorphism. Since $Ra\oplus J\cong R\cong Ra\oplus I$ 
cancellation (using Jordon-H\"older or Krull-Schmidt) shows that $J\cong I$. 
If $h:R\cong J\cong I$, then $R\stk{(f,g)}\lra Ra\oplus I\cong R$ is an \
isomorphism sending $1$ to $(a,i)$ to $a+i$, where $i=h(1)$. Hence it follows 
that $a+i$ is a unit. \hb \vp 
\begin{lm} \label{swan4} Let $R$ be a semisimple artinian ring and $I$ be a 
left ideal of $R$. Let $J=Ra+I$. Write $J=Re$, where $e$ is an idempotent $($possible since $J$ is projective. For detail cf. \cite{BK} Theorem 4.2.7$)$. 
Then there is an element $i\in I$ such that $a+i=ue$, where $u$ is a unit 
in $R$.  
\end{lm} 
{\bf Proof.} Since $R=J+R(1-e)=Ra+I+R(1-e)$, using Lemma \ref{B} 
we can find a unit $u=a+i+x(1-e)$ in $R$ for some $x\in R$. Since $a+i\in Re$, it follows that 
$ue=a+i$. \hb  

\begin{co} \label{swan5} 
Let $R$ be a semisimple artinian ring and $(a_1,\dots,a_n)^t$ be a 
column vector over $R$, where $n\ge 2$. Let $\Sigma Ra_i=Re$, where $e$ is 
an idempotent. Then there exists $\eps\in {\E}_n(R)$ such that 
$\eps (a_1,\dots,a_n)^t=(0,\dots,0,e)^t$. 
\end{co} 
{\bf Proof.} By Lemma \ref{swan4} we can write 
$ue=\Sigma_{i=1}^{n-1} b_ia_i+a_n$, 
where $u$ is a unit. Therefore, applying an elementary transformation we can 
assume that $a_n=ue$. Multiplying from the left by 
$(I_{n-2} \perp u\perp u^{-1})$ we can make $a_n=e$. Since all $a_i$ are 
left multiple of $e$, further elementary transformations reduce our vector 
to the required form. \hb \vp 

The following observation will be needed to do the case $2n=4$. 

\begin{lm} \label{swan6} Let $R$ be a semisimple artinian ring and $e$ be 
an idempotent. Let $f=1-e$, and $b$ be an element of $R$. If 
$fRb\subseteq Re$, then we have $b\in Re$. 
\end{lm}
{\bf Proof.} Since $R$ is a product of simple rings, it will suffice to do 
the case in which $R$ is simple. If $e=1$, we are done. Otherwise $RfR$ is a 
non-zero two sided ideal, and hence $RfR=R$. Since $Rb=RfRb\subseteq Re$, 
we have $b\in Re$. \hb 

\begin{lm} \label{swan7} Let $R$ be a semisimple artinian ring and let
$-:R\ra R$ be a $\lmd$-involution on $R$. Let 
$( x \,\, y)^t $ be a unimodular 
element of $R^{2n}$, where $2n\ge 4$. Then there exists an element $\eps\in 
{\E}(2n,R)$ such that $\eps( x \,\, y)^t=
(x' \,\, y')^t$, where $x_1'$ is a unit
in $R$.
\end{lm} 
{\bf Proof.} Let $x=(x_1,\ldots,x_r)^t$ and 
$b=(y_1,\ldots,y_r)^t$. We claim that there exists $\eps\in {\E}(2n,R)$ 
such that $\eps( x \,\, y)^t =
(x' \,\, y')$, where $x'$ is a unit in $R$.
Among all $( x' \,\, y')^t$ of this form, choose one
for which the ideal $I=\Sigma Rx_i'$ is maximal. Replacing the original 
$(x \,\, y )^t$ by $( x' \,\, y')^t$ we can assume that $I=\Sigma Rx_i$
is maximal among such ideals. Write $I=Re$, where $e$ is an idempotent in
$R$. By Corollary \ref{swan5} we can find an element $\eta\in {\E}_n(R)$ 
such that $\eta x =(0,0,\dots,e)^t$. Hence we assume that $x=(0,0,\dots,e)^t$.
We claim that $y_i\in Re$ for all $i\ge 1$.

First we consider the case $2n\ge 6$. 
Assume $y_1\ne I$, but $y_i\in I$ for all $i\ge 2$. If 
we apply $q\eps_{1n}(1)$ in the quadratic case
then this replaces $y_{n}$ to $y_{n}-y_1$ but 
not changes $e$ and $y_1$. On the other hand for the hermitian case we do not
have the generator $q\eps_{1n}(1)$. But if we apply $hm_n(1,\ldots,1)$, then
it changes $y_2$ but does not changes $e$ and $b_1$. 
Therefore, in both the cases we can therefore assume that some 
$y_i$ with $i>1$ is not in $I$. (Here recall that we have put no restriction 
on $C$, {\it i.e.} for us $C=R^r$). 
Apply $qr_{i i}(1)$ with $2\le i\le n$ in the quadratic case. 
This changes $x_i=0$ (for $i>1$) to $y_i$ while $x_n=e$ is
preserved. The ideal generated by the entries of $x$ now contains $Re+Ry_i$, 
which is larger than $I$, a contradiction, as $I$ is maximal. 
In the hermitian case if we apply suitable $hr_i(1,\ldots,1)$ then also we 
see that the ideal generated by the entries of $x$ now contains $Re+Ry_i$, 
hence a contradiction. 

If $2n=4$, we can argue as follows. Let $f=1-e$. Let us assume that 
$y_1\ne I$ as above. Then by Lemma \ref{swan6} it will follow that we can find
some $s\in R$ such that $fsy_1\ne Re$. 
First consider the quadratic case. 
Applying $qr_{21}(fs)$ replaces $x_2=e$ 
by $c=e+fsb_1$. As $ec=e$, $I=Re\subset Rc$. Also, $fc=fsb_1\in Rc$ but 
$fc\notin I$. Hence $I \subsetneq Rc$, a contradiction. We can get the similar 
contradiction for $y_2$  by applying $qr_{22}(fs)$. In the hermitian case,
apply $hr_1(1)$ to get the contradiction for $y_1$. Now
note that in this $r=1$ as we have assume $r<n$. Hence 
we can apply $qr_{22}(fs)$ to get the contradiction.

Since all $y_i$ lie in $Re$, the left ideal generated by the all entries of 
$(x\,\,y)^t$ is $Re$, but as this column 
vector is unimodular $Re=R$, and therefore $e=1$. \hb \vp \\
{\bf Proof of Theorem \ref{swan}.}  Let $J$ be the Jacobson radical of
$R$. Since the left and the right Jacobson radical are same, $J$ is stable
under the involution which therefore passes to $R/J$. Let $\eps$ be as in
Lemma \ref{swan7} for the image 
$(x'\,\,y')^t$ of
$(x\,\,y)^t$. By lifting $\eps$ from 
$R/J$ to $R$ and applying it to 
$(x\,\,y)^t$ we reduce to the case where 
$x_n$ is a unit in $R$. Let $\alpha=x_n\perp x_n^{-1}$. Then applying 
$(I_{n-2}\perp  \alpha \perp I_{n-2} \perp \alpha^{-1})$ we can assume that 
$x_n=1$.

Next applying $\Pi_{i=1}^{n-1} ql_{ni}(-y_i)$ and  
$\Pi_{i=1}^{n-1} hl_{ni}(-y_i)$ in the respective cases we
get $y_1=\cdots=y_{n-1}=0$. As 
isotropic vector remains isotropic under
elementary quadratic (hermitian) transformation, we have 
$y_n+\lambda\ol{y}_n=0$, hence $ql_{11}(\ol{\lmd}\ol{y}_n)$ and 
$hl_{11}(\ol{\lmd}\ol{y}_n)$ are defined and 
applying it reduces $y_n$ to $0$ in both the cases. 
Now we want to make $x_i=0$ for $i=1,\ldots,n$. In the quadratic case it can
be done by applying $\Pi_{i=1}^{n-1} h\eps_{in}(-x_i)$. Note
that this transformation does not affect any $y_i$'s, as $y_i=0$. In the
hermitian case we can make $x_{r+1}=\cdots=x_n=0$ as before applying 
$\Pi_{i=r+1}^{n-1} q\eps_{in}(-x_i)$. To make 
$x_1=\cdots=x_r=0$ we have to recall that the set $C=R^r$, {\it i.e.} there is
no restriction on the set $C$. Hence $hr_n(-x_1,\ldots,-x_r)$ is defined and
applying it we get $x_1=\cdots=x_r=0$. Also note that other $x_i$'s and
$y_i$'s remain unchanged. Finally, applying $hl_{nn}(1)$ and then
$hr_{nn}(-1)$ we get the required vector 
$(0,\ldots,0,1)$. This completes the proof.\hb 

\begin{tr}  
Let $R$ be a commutative ring with identity
and $A$  an associative $R$-algebra such that $A$ is finite 
as a left $R$-module. Then the following are equivalent for $n\ge 3$ in the quadratic case and $n\ge r+3$ in the hermitian case:
\begin{enumerate}
 \item {\bf (Normality)} ${\E}(2n, A, \LMD)$ is a normal subgroup of 
${\G}(2n, A, \LMD)$.
\item {\bf (L-G Principle)}
If $\alpha(X)\in {\G}(2n,A[X], \LMD[X])$, $\alpha(0)={\I}_n$ and 
{\tiny $$\alpha_{\m}(X)\in {\E}(n,A_{\m}[X], \LMD_{\mf m}[X])$$} for every maximal ideal $\m \in \M(R)$, then {\tiny $$\alpha(X)\in {\E}(2n,A[X], \LMD[X]).$$} $($Note that $A_{\m}$ denotes $S^{-1}A$, where 
$S = R - \m$.$)$
\end{enumerate}
\end{tr} 
{\bf Proof.}  
We have proved the Lemma \ref{key5} for any form ring with identity. In particular, suppose ${\E}(2n, A, \LMD)$ is a normal subgroup of ${\G}(2n, A, \LMD)$. 
Let $\alpha=I_n+M(v,w)$, where $v=Ae_1$, and $A\in {\G}(n,A,\LMD)$. 
Then we can write  $\alpha=A(I_n+M(e_1, w_1))A^{-1}$, where 
$w_1= A^{-1} w$.
Hence it is enough to show that 
$I_n + M(e_1,w_1)\in {\E}(n,R)$. Now arguing as in the proof of Lemma 
\ref{key5} we get the result. Now in section \S we have proved Local-Global Principle as a consequence of Lemma \ref{key5}. Hence the implication follows. 

To prove the converse we need $A$ to be finite as $R$-module, where $R$ is a commutative ring with identity ({\it i.e.} a ring with trivial involution). 

Let $\alpha\in {\E}(2n,A,\LMD)$ and $\beta\in {\G}(2n,A, \LMD)$. Then 
$\alpha=\Pi ge_{ij}(x)$, $x\in A$. Hence, 
$\beta\alpha \beta^{-1}= \Pi (I_{2n}+x\beta\, M(\star_1,\star_2)\beta^{-1})$, where $\star_1$ and $\star_2$ are suitably chosen standard basis vectors. Now let $v=\beta e_i$ and $w = x \beta  e_j$. 
Then we get $\beta\alpha \beta^{-1}= \Pi (I_{2n}+M(v,w))$,  
where $v\in {\rm Um}_{2n}(A)$ and $\langle v,w\rangle=0$. 
We show that each $(I_{2n}+M(v,w))\in {\E}(2n,A,\LMD)$.

Let $\gamma(X)=I_{2n}+XM(v,w)$. Then $\gamma(0)=I_{2n}$. 
By Lemma \ref{HB} it follows that $S^{-1}A$ is a semilocal ring, where $S=R-\m$, 
$\m \in {\M}(R)$. Since $v\in {\rm Um}_{2n}(A)$, using Theorem \ref{swan} we get 
$v\in {\E}(2n, S^{-1}A, S^{-1}\LMD) e_1$, hence $Xv\in {\E}(2n, S^{-1}A[X], S^{-1}\LMD[X]) e_1$. Therefore, applying Lemma \ref{key5} over 
$S^{-1}(A[X], \LMD[X])$ it follows that {\tiny $$\gamma_{\m}(X)\in {\E}(2n,S^{-1}A[X],S^{-1}\LMD[X]).$$} Now applying Theorem \ref{LG} it follows that 
$\gamma(X)\in {\E}(2n,A[X],\LMD[X])$.
Finally, putting $X=1$ we get the result. \hb

\section{Nipotent property for ${\k}$ of Hermitian groups}

We devote this section to discuss the study of nilpotent property of unstable ${\k}$-groups. 
The literature in this direction can be seen in the work of 
A. Bak, N. Vavilov and R. Hazarat and that we have already discussed in the Introduction. 
Throughout this section we assume $R$ is a commutative ring with identity, {\it i.e.} we are considering trivial involution and $n\ge r+3$. Following is the statement of the theorem. 

\begin{tr} \label{nil}
The quotient 
group $\frac{{\SH}_{2n}(R, a_1,\ldots, a_r)}{{\EH}_{2n}(R, a_1,\ldots, a_r)}$ is nilpotent for $n\ge r+3$. The  class of nilpotency is at the most \tn{max}  $(1, d+3-n)$,  where $d=\dim \,(R)$.
\end{tr}

The proof follows by emitting the proof of Theorem 4.1 in \cite{BBR}. 

\begin{lm} \label{sol3a}
Let $\beta\in {\SH}(2n,R, \LMD)$, with $\beta\equiv I_n$ modulo $I$, 
where $I$ is an ideal contained in the Jacobson 
radical $J(R)$ of $R$. Then there 
exists $\theta\in {\EH}_{2n}(R, a_1,\ldots, a_r)$ such that $\beta\theta$= the diagonal matrix 
$[d_1,d_2,\dots,d_{2n}]$, where each $d_i$ is a unit in $R$ with $d_i\equiv 1$ 
modulo $I$, and $\theta$ a product of 
elementary generators with each congruent to identity modulo $I$.
\end{lm} 
{\bf Proof.} The diagonal elements are units. Let $\beta=(\beta_{ij})$, where 
$d_i=\beta_{ii}=1+s_{ii}$ with $s_{ii}\in I\subset J(R)$ for $i=1,\ldots, 2n$, and $\beta_{ij}\in I\subset J(R)$ for $i\ne j$. First we make all the $(2n,j)$-th, and $(i,2n)$-th  entries zero for $i=2,\ldots,n$, $j=2,\ldots,n$. 
Then repeating the above process we can reduce the size of $\beta$. Since we are considering trivial involution, we take {\tiny $$\alpha
=\underset{j=1}{\overset{n}\Pi}hl_{nj}(-\beta_{2nj}d_{j}^{-1})
\underset{n+1\le j\le n+r}{\underset{n+r+1\le i\le 2n-1}\Pi}hm_i(-\zeta_jd_j^{-1})
\underset{n+r+1\le j\le 2n-1}{\underset{r+1\le i\le n-1} \Pi} h\eps_{in}(\beta_{\rho(n)\rho(i)}d_{j}^{-1}),$$}
where $j=i-r$ and $\zeta_j=(0,\ldots,0,\beta_{2n  j})$, and \\
{\tiny $$\gamma=\underset{r+1\le i\le 2n-1}{\underset{r+1\le j\le 2n-1}\Pi h\eps_{nj}(a_{i-r}(\star)d_{2n}^{-1})}
hr_n(\eta),$$} where $a_t=0$ for $t>r$, and 
$\eta=(\beta_{1 2n}d_{2n}^{-1},\beta_{2 2n}d_{2n}^{-1},\ldots,\beta_{n 2n}d_{2n}^{-1})$.
Then the last column and last row of $\gamma\beta\alpha$ become 
$(0,\dots,0, d_{2n})^t$, where $d_{2n}$ is a unit in $R$ and $d_{2n}\equiv 1$ modulo 
$I$. Repeating the process we can modify $\beta$ to the required form. \hb

\begin{lm} \label{sol3}
Let $(R, \LMD)$ be a commutative form ring, {\it i.e.} with trivial involution and $s$ be a non-nilpotent element in $R$.
Let $D$ denote the diagonal matrix $[d_1,\dots,d_{2n}]$, where 
$d_i\equiv 1$ modulo $(s^l)$ for $l\ge 2$. Then 
{\tiny $$\left[ge_{ij}\left(\frac{a}{s} X \right), D\right]\subset 
{\EH}_{2n}(R[X], a_1,\ldots, a_r)\cap {\SH}_{2n}((s^{l-1})R,  a_1,\ldots, a_r).$$}
\end{lm} 
{\bf Proof.} Let $d=d_id_j^{-1}$. Then using a 
list of commutator laws for elementary generators is stated in (\cite{T}, pg. 237-239) for the hermitian group, it follows that  
{\tiny $$\left[ge_{ij}\left(\frac{a}{s} X \right), D\right]= 
ge_{ij}\left(\frac{a}{s} X\right) ge_{ij}\left(-\frac{a}{s} dX\right).$$} Since 
$d_i, d_j\equiv 1$ modulo $(s^l)$ for $l\ge 2$, we can write $d=1+s^m\lambda$ 
for some $m>2$ and $\lambda\in R$. Hence 
{\tiny \begin{align*}
ge_{ij}\left(\frac{a}{s} X \right)ge_{ij}\left(-\frac{a}{s} dX\right) & =
ge_{ij}\left(\frac{a}{s} X\right) ge_{ij}\left(-\frac{a}{s} X\right)
ge_{ij}\left(-\frac{a}{s}s^m\lambda X\right)\\
=  ge_{ij}\left(-\frac{a}{s}s^m\lambda X\right)\in &
{\EH}_{2n}(R[X], a_1,\ldots, a_r)\cap {\SH}_{2n}((s^{m-1})R, a_1,\ldots, a_r).
\end{align*}} \hb 
\begin{lm} \label{sol4}
Let $(R, \LMD)$ be as above, $s\in R$ a non-nilpotent element in $R$ and $a\in R$. 
Then for $l\ge 2$  
{\tiny $$\left[ge_{ij}\left(\frac{a}{s} X\right), {\SH}_{2n}(s^lR, a_1,\ldots, a_r)\right]
\subset {\EH}_{2n}(R[X], a_1,\ldots, a_r).$$}
More generally, 
$\left[\eps(X), {\SH}_{2n}(s^lR[X], a_1,\ldots, a_r)\right]\!\!\subset \!\!{\EH}_{2n}(R[X],\!\! a_1,\ldots,\!\! a_r )$ for $l\gg 0$ and 
$\eps(X)\in {\EH}_{2n}(R_s[X], a_1,\ldots, a_r)$. 
\end{lm} 
{\bf Proof.} First fix $(i,j)$ for $i\ne j$.
 
Let $\alpha(X)=[g_{ij}\left(\frac{a}{s} X\right), \beta]$ for 
some $\beta\in {\SH}_{2n}(s^lR, a_1,\ldots, a_r)$. 
Since $s$ is in $J(R)$, the diagonal entries of $\beta$ are unipotent.
As $l\ge 2$, it follows that 
$\alpha(X)\in {\SH}_{2n}(R[X], a_1,\ldots, a_r)$. Since ${\EH}_{2n}(R[X], a_1,\ldots, a_r)$ is a normal subgroup of 
${\SH}_{2n}(R[X], a_1,\ldots, a_r)$ for $n\ge r+3$ ({\it cf.}\cite{T}, Theorem 4.2), we get {\tiny $$\alpha_s(X)\in {\EH}_{2n}(R_s[X], a_1,\ldots, a_r).$$ }
Let  $B=1+sR$. We show that $\alpha_B(X)\in {\EH}_{2n}(R_B[X], a_1,\ldots, a_r)$. Since 
$s\in J(R_B, \LMD_B)$,  it follows from Lemma \ref{sol3a}
that we can decompose $\beta_B=\eps_1\cdots\eps_tD,$ where 
$\eps_i=ge_{p_iq_i}(s^l \lambda_i) \in {\EH}_{2n}(R_B, a_1,\ldots, a_r)$; $\lambda_i\in R_B$
and $D$ = the diagonal matrix $[d_1,\dots,d_{2n}]$ with $d_i$ is a unit in $R$ 
and $d_i\equiv 1$ modulo $(s^l)$ for $l\ge 2$; $i=1,\dots,2n$. 
If $t=1$, then using the commutator laws for elementary generators is stated in (\cite{T}, pg. 237-239)
it follows from Lemma \ref{sol3}
that $\alpha_B(X)\in {\EH}_{2n}(R_B[X], a_1,\ldots, a_r)$. Suppose $t>1$. Then 
{\tiny \begin{align*}
\alpha_B(X) &
 =\left[ge_{ij}\left(\frac{a}{s} X\right), \eps_1 \cdots \eps_t D\right] \\
&  =\left[ge_{ij}\left(\frac{a}{s} X\right), \eps_1\right] \eps_1 
\left[ge_{ij}\left(\frac{a}{s} X\right),\eps_2 \cdots \eps_t D\right] 
\eps_1^{-1}
\end{align*}}
and by induction each term is in ${\EH}_{2n}(R_B[X], a_1,\ldots, a_r)$, hence 
{\tiny $$\alpha_B(X)\in {\EH}_{2n}(R_B[X], a_1,\ldots, a_r).$$} Since $\alpha(0)=I_n$, 
by the Local-Global Principle for the hermitian groups (Theorem \ref{LG}) it follows that 
$\alpha(X)\in {\EH}_{2n}(R[X], a_1,\ldots, a_r)$. \hb
\begin{co} \label{sol5}
Let $R$ be as above, $s\in R$ be a non-nilpotent element in $R$ and $a\in R$. 
Then for $l\ge 2$  
{\tiny $$\left[ge_{ij}\left(\frac{a}{s} \right), {\SH}_{2n}(s^lR, a_1,\ldots, a_r)\right]
\subset {\EH}_{2n}(R, a_1,\ldots, a_r).$$}
More generally, 
$\left[\eps, {\SH}_{2n}(s^lR, a_1,\ldots, a_r)\right]\subset {\EH}_{2n}(R, a_1,\ldots, a_r)$ for $l\gg 0$ and $\eps\in {\EH}_{2n}(R_s, a_1,\ldots, a_r)$. 
\end{co}  
{\bf Proof of Theorem \ref{nil}:}
Recall 

 Let $G$ be a group. Define $Z^0=H$, $Z^1=[G,G]$ and 
$Z^i=[G,Z^{i-1}]$. Then $G$ is said to be nilpotent if $Z^r=\{e\}$ for
some $r>0$, where $e$ denotes the identity element of $H$.

Since the map ${\EH}_{2n}(R, a_1,\ldots, a_r)\ra {\EH}_{2n}(R/I, \ol{a}_1,\ldots, \ol{a}_r)$
is surjective we may and do assume that $R$ is a reduced 
ring. Note that if $t\ge d+3$, then the group ${\SH}_{2n}(R, a_1,\ldots, a_r)/{\EH}_{2n}(R, a_1,\ldots, a_r)={\KH}_1(R, a_1,\ldots, a_r)$, 
which is abelian and hence nilpotent. So we consider the case $t\le d+3$. 
Let us first fix a $t$. We prove the theorem by induction 
on $d=\dim R$. Let $G={\SH}_{2n}(R, a_1,\ldots, a_r)/{\EH}_{2n}(R, a_1,\ldots, a_r)$. Let $m=d+3-t$ and 
$\alpha=[\beta, \gamma]$ for some $\beta\in G$ and $\gamma\in Z^{m-1}$. 
Clearly, the result is true for $d=0$. 
Let $\widetilde{\beta}$ be the pre-image of $\beta$ under the map 
${\SH}_{2n}(R, a_1,\ldots, a_r)\ra {\SH}_{2n}(R, a_1,\ldots, a_r)/{\EH}_{2n}(R, a_1,\ldots, a_r)$.
If $R$ is local then arguing as Lemma \ref{sol3a} is follows that $\EH_{2n}=\SH_{2n}$, hence we can 
choose a non-zero-divisor $s$ in $R$ such that  
$\widetilde{\beta}_s\in {\EH}_{2n}(R_s, a_1,\ldots, a_r)$.

Consider $\ol{G}$,
where bar denote reduction modulo $s^l$ for some $l\gg 0$. By the  
induction hypothesis $\ol{\gamma}=\{1\}$ in $\ol{\GH}_{2n}$. 
Since ${\EH}_{2n}$ is a normal subgroup
of ${\SH}_{2n}$ for $n\ge r+3$, 
by modifying $\gamma$ we may assume that
$\widetilde{\gamma}\in {\SH}_{2n}(R,s^lR, a_1,\ldots, a_r)$, where $\widetilde{\gamma}$ is the pre 
image of $\gamma$ in ${\SH}_{2n}(R,a_1,\ldots, a_r)$. Now by  Corollary \ref{sol5} it follows that 
$[\widetilde{\beta},\widetilde{\gamma}]\in {\EH}_{2n}(R,a_1,\ldots, a_r)$. 
Hence $\alpha=\{1\}$ in $G$. \hb 

{\tiny \begin{re} \tn{In (\cite{BRK}, Theorem 3.1) it has been proved that 
the question of normality of the elementary subgroup and the Local-Global Principle are 
equivalent for the elementary subgroups of the linear, symplectic and 
orthogonal groups over an almost commutative ring with identity. 
There is a gap in the proof of the statement $(3)\Ra (2)$ of Theorem 3.1 in 
\cite{BRK} (for an almost commutative ring). 
The fact that over a non-commutative semilocal ring 
the elementary subgroups of the classical groups 
acts transitively on the set of unimodular and istropic 
({\it i.e.} $\langle\,v,v\rangle=0$) vectors of length $n\ge 2$ in the linear
case and $n=2r\ge 4$ in the non-linear cases has been used 
in the proof, but it is not mentioned anywhere in the article. 
This was pointed by Professor R.G. Swan and he provided us a proof for the above result. }
\end{re}
 {\bf Acknowledgment:} My sincere thanks to Professors R.G. Swan for giving me his permission to
reproduce his proof of Theorem \ref{swan} (he gave a proof for the symplectic and orthogonal groups as noted above). (The proof of H. Bass' Lemma \ref{B} given here is a simple proof of it that he used to give in his K-theory courses). I thank 
Professors A. Bak and Gouping Tang for their important remarks, and Professor Ravi A. Rao for his continuous encouragement. I have started this work as a Post-doctoral fellow at Harish Chandra Research Institute, Allahabad. I really appreciate the ambiance and facilities at HRI. 
I am grateful to Stat-Math Unit, Indian Statistical Institute, Kolkata for allowing me to use their infrastructure facilities even after my Post-doctoral period. My special thanks to Professor S.M. Srivastava for his concern regarding my research. Finally, I thank Gaurab Tripathi for helping me to correct the final manuscript. }

{\tiny
 
\addcontentsline{toc}{chapter}{Bibliography} 
{\scshape Indian Institute of Science Education and Research - Kolkata},\\
{\scshape Mohanpur Campus, P.O. -- BCKV Campus Main Office},\\
{\scshape Mohanpur -- 741252, Nadia,West Bengal, India}.\\
{\it Email: rabeya.basu@gmail.com, rbasu@iiserkol.ac.in}}

\end{document}